\documentclass[a4paper]{amsart} 

\usepackage{amsmath, amscd, amsthm, amssymb, mathrsfs}

\DeclareMathOperator{\Aut}{Aut}

\DeclareMathOperator{\tr}{tr}
\DeclareMathOperator{\rank}{rank}

\DeclareMathOperator{\Scal}{Scal}

\DeclareMathOperator{\CM}{CM}

\newcommand{\del}{\partial}
\newcommand{\delb}{\bar{\del}}

\theoremstyle{plain}
	\newtheorem{theorem}{Theorem}[section]
	
	\newtheorem{lemma}[theorem]{Lemma}
	\newtheorem{corollary}[theorem]{Corollary}

\theoremstyle{definition}

\theoremstyle{plain}
	\newtheorem*{theorem*}{Theorem}
	\newtheorem*{proposition*}{Proposition}
	\newtheorem*{lemma*}{Lemma}
	\newtheorem*{corollary*}{Corollary}
	\newtheorem*{conjecture*}{Conjecture}
\theoremstyle{definition}
	\newtheorem*{definition*}{Definition}
	\newtheorem*{remark*}{Remark}
	\newtheorem*{remarks*}{Remarks}

\numberwithin{equation}{section}

\begin{document}

\title[Fibrations with cscK metrics]{Fibrations with constant scalar curvature K\"ahler metrics and the CM-line bundle}
\author{Joel Fine}

\begin{abstract}
Let $\pi \colon X \to B$ be a holomorphic submersion between compact K\"ahler manifolds of any dimensions, whose fibres and base have no non-zero holomorphic vector fields and whose fibres admit constant scalar curvature K\"ahler metrics. This article gives a sufficient topological condition for the existence of a constant scalar curvature K\"ahler metric on $X$. The condition involves the $\CM$-line bundle---a certain natural line bundle on $B$---which is proved to be nef. Knowing this, the condition is then implied by  $c_1(B) <0$. This provides infinitely many K\"ahler manifolds of constant scalar curvature in every dimension, each with K\"ahler class arbitrarily far from the canonical class.
\end{abstract}

\maketitle

\section{Introduction}

\subsection{Background}

In \cite{calabi:ekm} Calabi proposed that, when one exists, a constant scalar curvature K\"ahler (cscK) metric should provide a canonical representative for a given K\"ahler class. Since this suggestion, much work has focused on the topic. The general existence theory has been looked at in depth, motivated by a suggestion of Yau \cite{yau:opig} relating the existence of a K\"ahler--Einstein metric on a Fano manifold to the stability of the canonical polarisation in the sense of geometric invariant theory. Later, Tian \cite{tian:kemwpsc} introduced the notion of K-stability, and suggested that the existence of a cscK metric in the first Chern class of an ample line bundle $L$ should be equivalent to the K-stability of the polarisation defined by $L$. Finally, Donaldson \cite{donaldson:fml} gave an interpretation of the conjecture in terms of symplectic geometry and moment maps. See \cite{donaldson.scasotv} for a precise statment of the conjecture, including a definition of K-stability. 

The difficulty, from the analytic viewpoint, in determining whether or not a cscK metric exists is that the resulting PDE is fourth order and fully non-linear. Despite this, some examples have been found, but only in situations with prescribed geometry which is then used get a handle on the PDE. For two recent such examples see \cite{apostolov:h2fikg3ce,rollin.singer:cokswcsc}. 

This paper adds to the list of examples, the prescribed geometry here being that of a fibration. CscK metrics are found on  the total space of a holomorphic submersion $X \to B$ between compact K\"ahler manifolds. Even in this restricted context, however, a clear relationship between the existence of a cscK metric and stability arises. Moreover, the varying moduli of the fibres plays an important r\^ole. One is led to a natural line bundle on $B$---the so-called $\CM$-line bundle---which is shown to be nef on the moduli space of K-semi-stable polarised varieties (ignoring the technical issues of whether such a space exists).

\subsection{Overview of the results}

Let $\pi \colon X \to B$ be a holomorphic submersion between compact K\"ahler manifolds and let $L\to X$ be a relatively ample line bundle. There is a natural line bundle, called the CM-line bundle, on $B$ which was introduced by Tian \cite{tian:tkeohas}. It is defined as follows.

Let $Y$ be a fibre, $n = \dim Y$ and let $s$ denote the constant
\begin{equation}\label{s}
s = \frac{n c_1(Y) c_1(L)|_Y^{n-1}}{c_1(L)|_Y^n}.
\end{equation}
(In terms of differential geometry, given a K\"ahler metric in $c_1(L)$, its restriction to a fibre has average scalar curvature $2\pi s$. In terms of algebraic geometry, $s$ is essentially the second coefficient in the Hilbert polynomial of $L|_Y$.) Let $K_{X/B}$ denote the relative canonical bundle, and let $\mathcal E$ denote the virtual bundle
$$
\mathcal E
=
(n+1)(K_{X/B}^* - K_{X/B})\otimes(L - L^*)^n - s(L-L^*)^{n+1}.
$$
The CM-line bundle is then $L_{\CM} = [\det(\pi_!E)]^*$. The bundle makes sense when $\pi$ is a flat proper morphism of varieties. This level of generality is not needed for the result concerning existence of cscK metrics. 

In fact, only the first Chern class of $L_{\CM}$ is used in this article. It can be defined without recourse to virtual bundles or direct images: by Grothendieck--Riemann--Roch,
$$
c_1(L_{\CM}) =
2^{n+1}
\pi_*\left[ 
\left((n+1)c_1(K_{X/B}) + s c_1(L) \right) c_1(L)^n 
\right ].
$$
Notice that $c_1(L_{\CM})$ is unchanged when $L$ is twisted by a line bundle pulled back from $B$. This is because after pushing forward the factor involving $s$ cancels other unwanted terms.

The relation of $L_{\CM}$ to K-stability and the existence of cscK metrics has already been studied by Tian, \cite{tian:kemwpsc} and Paul--Tian \cite{paul.tian:aogs,paul.tian:aaaks}. The fact that $c_1(L_{\CM})$ is well behaved when all fibres of $X$ admit a cscK metric is noted in Fujiki--Schumacher \cite{fujiki.schumacher:tmsoeckmagwpm}.

For notational convenience, let $\alpha$ denote the following class in $H^{1,1}(B)$:
\begin{equation}\label{alpha} 
\alpha
=
\frac{c_1(L_{\CM})}{2^{n+1}(n+1)\pi_*(c_1(L^n))}.
\end{equation}
The main theorem proved here is the following. (See section \ref{proof of main theorem} for the proof.)

\begin{theorem}\label{main theorem}
Let $\pi\colon X\to B$ be a holomorphic submersion between compact K\"ahler manifolds whose fibres and base admit no non-zero holomorphic vector fields. Let $L$ be a relatively ample line bundle on $X$ such that the restriction of $c_1(L)$ to each fibre admits a cscK metric. Suppose, moreover, that $\alpha -c_1(B) \geq 0$.

Then, for all sufficiently large $r$, the class 
$$
\kappa_r = c_1(L) + r \pi^*\kappa_B
$$
contains a cscK metric where, if $\alpha -c_1(B)=0$, $\kappa_B$ is any K\"ahler class on the base, whilst if $\alpha -c_1(B) >0$, then $\kappa_B=\alpha -c_1(B)$.
\end{theorem}

In checking the condition $\alpha -c_1(B) \geq0$, it is helpful to know about the positivity of $\alpha$ itself. This is provided by the following result. (See section \ref{non-negativity of alpha} for the proof.) Here, a line bundle $E \to Y$ is \emph{asymptotically Hilbert stable} if for all sufficiently large $m$, $E^m$ gives a Hilbert stable embedding of $Y$.

\begin{theorem}\label{alpha is nef}
Let $\pi\colon X\to B$ be a flat proper morphism of varieties. Let $L$ be a relatively ample line bundle on $X$ whose restriction to a generic fibre is asymptotically Hilbert semi-stable. 

Then $L_{\CM}$ is nef; that is, $c_1(L_{\CM})$ evaluates non-negatively over any curve in $B$.
\end{theorem}

It is a result due to Donaldson \cite{donaldson:scape} that if $E \to Y$ is an ample line bundle over a compact K\"ahler manifold with $\Aut(Y,E)$ discrete (modulo scalars), then the existence of a cscsK metric in $c_1(E)$ implies that $E$ is asymptotically Hilbert stable (and hence semi-stable). Combining this with Theorems \ref{main theorem} and \ref{alpha is nef} and the fact that nef plus positive is positive gives the following corollary.

\begin{corollary}\label{anti-fano base}
Let $\pi \colon X \to B$ be a holomorphic submersion of compact K\"ahler manifolds whose fibres and base admit no non-zero holomorphic vector fields. Let $L$ be a relatively ample line bundle on $X$ such that the restriction of $c_1(L)$ to each fibre admits a cscK metric. Suppose, moreover, that $c_1(B) <0$. 

Then, for all sufficiently large $r$, the class
$$
\kappa_r = c_1(L) + r \pi^*(\alpha - c_1(B))
$$
contains a cscK metric. 
\end{corollary}

The condition $c_1(B) <0$ is straight forward to check and general enough to provide a large number of examples of cscK metrics. Some of these are described in section \ref{examples}, giving infinitely many K\"ahler manifolds of constant scalar curvature in each dimension, all with K\"ahler class arbitrarily far from the canonical class.

When $c_1(L_{\CM})$ is actually ample, one can say something in the case when $c_1(B)=0$. In this direction, Fujiki--Schumacher \cite{fujiki.schumacher:tmsoeckmagwpm} prove the following.

\begin{theorem}[Fujiki--Schumacher]\label{Fujiki--Schumacher}
Let $\pi \colon X \to B$ be a holomorphic submersion whose fibres admit no non-zero holomorphic vector fields. Let $L \to X$ be a relatively ample line bundle such that the restriction of $c_1(L)$ to each fibre admits a cscK metrics. Suppose, moreover, that $\pi$ is not trivial over any curve in $B$. Then $c_1(L_{\CM})$ is ample on $B$.
\end{theorem}

(Here $\pi$ is trivial over a curve means that all fibres of $\pi$ over that curve are biholomorhpic as polarised varieties.) Combining this with Theorem \ref{main theorem} gives the following.

\begin{corollary}
Let $\pi \colon X \to B$ be a holomorphic submersion of compact K\"ahler manifolds whose fibres and base admit no non-zero holomorphic vector fields. Let $L$ be a relatively ample line bundle on $X$ such that the restriction of $c_1(L)$ to each fibre admits a cscK metric. Suppose, moreover, that $\pi$ is not isotrivial on any curve in $B$ and that $c_1(B) =0$. 

Then, for all sufficiently large $r$, the class
$$
\kappa_r = c_1(L) + r \pi^*\alpha
$$
contains a cscK metric. 
\end{corollary}

Note that neither of Theorems \ref{alpha is nef} and \ref{Fujiki--Schumacher} imply the other; whilst ample is stronger than nef, Theorem \ref{alpha is nef} makes no mention of cscK metrics and applies both when the fibres are singular and have automorhpisms.

\subsection{Acknowledgements}

I would like to thank Simon Donaldson, Julius Ross and Richard Thomas for helpful conversations concerning this paper and related matters. I am also grateful to Bill Harvey for bringing to my attention the iterated surface bundles in section \ref{examples}.

\section{Proof of Theorem \ref{main theorem}}\label{proof of main theorem}

Assume throughout this section that $\pi \colon X \to B$ is a holomorphic submersion between compact K\"ahler manifolds, that the fibres and base of $\pi$ have no non-zero holomorphic vector fields, that $L$ is relatively ample and that the restriction of $c_1(L)$ to each fibre admits a cscK metric. Moreover, assume that $\alpha - c_1(B) \geq 0$. If $\alpha - c_1(B) = 0$, let $\kappa_B$ denote any K\"ahler class on $B$; if $\alpha -c_1(B) >0$, let $\kappa_B = \alpha - c_1(B)$.

For sufficiently large $r$, $\kappa_r = c_1(L) + r\pi^*\kappa_B$ is a K\"ahler class on $X$. The first step is to find a K\"ahler representative whose fibrewise restriction is cscK.  Begin with any K\"ahler form $\omega' \in \kappa_{r_0}$ for some large fixed $r_0$. 

For $b \in B$, write $Y_b = \pi^{-1}(b)$ and let $\omega_b$ be a cscK metric in $c_1(L)|_{Y_b}$. Since $H^0(TY_b) = 0$, a theorem of Donaldson \cite{donaldson:scape} says that $\omega_b$ is in fact unique. For each $b$, there is a unique function $\phi_b \in C^\infty(Y_b)$ with mean-value zero (with respect to $\omega_b$) and such that $\omega'|_{Y_b} + i \delb \del \phi_b = \omega_b$. Provided the $\phi_b$ are smooth in $b$, they fit together to give a function $\phi \in C^\infty(X)$; then $\omega''=\omega' + i \delb\del \phi$  is a $(1,1)$-form in $\kappa_{r_0}$ whose fibrewise restriction is cscK. The smoothness of $\omega''$ is provided by the following lemma.

\begin{lemma}
The functions $\phi_b$ depend smoothly on $b$. 
\end{lemma}

\begin{proof}
This is essentially a standard result in the theory of elliptic PDEs. 
Choose a local trivialisation for $\pi$ over some disc $D \subset B$, $X|_D \cong D \times Y$ as smooth manifolds. By restriction in this trivialisation, the complex structure and K\"ahler form $\omega'$ on $X$ give a smooth family of K\"ahler structures $(J_b, \omega'_b)$ on $Y$ parameterised by $b \in D$. 

Define a map $S \colon D \times C^\infty(Y) \to C^\infty(Y)$
by
$$
S(b,\phi) = \Scal(\omega'_b +i(\delb \del)_b \phi)
$$
where $(\delb \del)_b$ is defined with respect to $J_b$. $S$ extends to a smooth map $D \times L^2_{k+4}(Y) \to L^2_k(Y)$. ($S$ is smooth in the $D$ factor because $\Scal(\omega,J)$ depends smoothly on $\omega$ and $J$, see, for example, section 2.2 in \cite{fine:csckmofcs}.)

By definition, $S(b, \phi_b)$ is a constant. The linearisation of $S$ with respect to $\phi$ at such a point is given by$$
\mathcal D_b^* \mathcal D_b\colon C^\infty(Y) \to C^\infty(Y).
$$
Here $\mathcal D_b$ is defined by 
$$
\mathcal D_b (\psi) = \delb_b (h_b(\psi))
$$
where $h_b(\psi)$ is the Hamiltonian vector field of $\psi$ with respect to $\omega'_b$ and $\delb_b$ is the $\delb$-operator on $TY$ determined by $J_b$. The map $\mathcal D^*_b$ is the $L^2$-adjoint of $\mathcal D$ with respect to the metric $\omega'_b$. (The linearisation of the scalar curvature map is computed in several places in the literature; see, for example, \cite{fine:csckmofcs}.)

Now $\mathcal D^*_b\mathcal D_b$ is elliptic (again, see \cite{fine:csckmofcs}), self adjoint and so has index zero. Since the fibres of $\pi$ have no non-zero holomorphic vector fields, $\ker \mathcal D^*_b\mathcal D_b$ is the constant functions. Hence $\mathcal D^*_b\mathcal D_b$ is an isomorphism modulo the constants. By the implicit function theorem, the map $b \mapsto \phi_b$ is a smooth map $D \to L^2_k(Y)$ for any $k$. By Sobolev embedding, it is a smooth map $D \to C^r(Y)$ for any $r$. Hence $\phi_b$ is smooth in $b$.
\end{proof}

Now $\omega'' = \omega' + i\delb \del \phi$ needn't be K\"ahler as it may be degenerate transverse to the fibres. To fix this, take a K\"ahler form $\omega_B \in\kappa_B$; for sufficiently large $r$, the form $\omega=\omega'' + (r-r_0)\pi^*\omega_B$ is a K\"ahler metric in $\kappa_r$ whose fibrewise restriction is cscK. 

The idea is that for large $r$, the geometry is dominated by that of the fibres which are cscK. One might hope then to be able to perturb the metric $\omega$ to a genuine cscK metric. A previous paper \cite{fine:csckmofcs} considers this problem in detail. In the case where $X$ is a surface it solves it completely; in higher dimensions it reduces it to solving a certain PDE for a metric on $B$. To describe this PDE, some notation is required. 

The fibrewise cscK metrics define a Hermitian structure in the vertical tangent bundle $V$ and hence in the relative canonical bundle $K_{X/B} = \Lambda^{\text{max}}V^*$; denote its curvature by $F$. Notice that the restriction to a fibre of $-iF$ is the Ricci form of that fibre with its cscK metric. 

The metric $\omega$ defines a vertical-horizontal decomposition of $TX$. Let $F_H$ denote the horizontal-horizontal component of $F$ with respect to this splitting. Define a form $a \in \Omega^{1,1}(B)$ by taking the fibrewise mean value of $iF_H$ with respect to $\omega$. More precisely, the differential $D\pi$ identifies the horizontal distribution restricted to $Y_b$ with $Y_b \times T_bB$; using this identification, on $Y_b$, $F_H \colon Y_b \to \Lambda^{1,1}T^*_bB$ is a vector valued function; then $a_b \in \Lambda^{1,1}T^*_bB$ is given by
$$
a_b = 
\frac{\int_{Y_b} iF_H \omega_b^n}
{\int_{Y_b}\omega_b^n}.
$$

\begin{theorem}[\cite{fine:csckmofcs}]\label{old theorem}
Let $\pi\colon X\to B$ be a fibration and $L \to X$ a relatively ample line bundle as above. Suppose, moreover, that there is a K\"ahler metric $\omega_B$ on the base solving
\begin{equation}\label{traced equation}
\Scal(\omega_B) - \tr_{\omega_B}a = \lambda,
\end{equation}
for some constant $\lambda$, and that there are no other cohomologous solutions to this equation. Then, for all sufficiently large $r$, the class
$
\kappa_r = \kappa + r\pi^*[\omega_B]
$
contains a cscK metric.
\end{theorem}

The difficulty with applying Theorem \ref{old theorem} is that the PDE for $\omega_B$ is as awkward to solve as the cscK equation. However, the \emph{un-traced} version, 
\begin{equation}\label{un-traced equation}
\rho(\omega_B) - a = \lambda \omega_B,
\end{equation}
(where $\rho(\omega)$ denotes the Ricci form of $\omega$) is essentially the complex Monge--Ampere equation. When $\lambda \leq 0$ this has been solved by Aubin \cite{aubin:edtmaslvkc} and Yau \cite{yau:otrcoackmatcmae1}. Before describing this, it is first necessary to give a cleaner description of $a$. (See (\ref{s}) and (\ref{alpha}) for the definitions of $s$ and $\alpha$.)

\begin{lemma}
$$
a
=
\frac{1}
{\pi_*(\omega^n)}\,
\pi_*\left[ 
\left( i F + \frac{2\pi s}{n+1} \omega\right)\wedge \omega^n\right]
$$
In particular, $a$ is closed and $[a] = 2\pi \alpha$.
\end{lemma}

\begin{proof}
Write $F_{V}$ and $F_H$ for the purely vertical and purely horizontal components of $F$ respectively. Since $-iF_V$ is the Ricci-form of $\omega_b$ which is cscK, $-niF_V\wedge\omega_b^{n-1} = 2 \pi s \omega_b^n$. So,
\begin{eqnarray*}
 \pi_*(iF_V \wedge \omega^n)
&=&
\pi_*(niF_V\wedge\omega_b^{n-1}\wedge\omega_H)\\
&=&
-2 \pi s\pi_*(\omega_b^n\wedge\omega_H)\\
&=&
-\frac{2\pi s}{n+1}\pi_*(\omega^{n+1})
\end{eqnarray*}
Hence,
\begin{eqnarray*}
\pi_*(iF\wedge\omega^n) 
&=& 
\pi_*(iF_H\wedge \omega^n) + \pi_*(iF_V\wedge\omega^n)\\
&=&
\pi_*(\omega^n)a
-
\frac{2\pi s}{n+1}\pi_*(\omega^{n+1})
\end{eqnarray*}
which gives the formula for $a$. 

To deduce that $[a] = 2\pi\alpha$, observe that $[i F] = 2\pi c_1(K_{X/B})$, whilst $[\omega] = \kappa_r = c_1(L) + \pi^*\kappa_B$. Twisting $L$ by bundles pulled back from $B$ leaves $\alpha$ unchanged, so $\alpha$ can be computed using $\kappa_r$ in place of $c_1(L)$.
\end{proof}

Everything is now in place to complete the proof.

\begin{proof}[Proof of Theorem \ref{main theorem}]

Let $\omega'_B \in \kappa_B$ be a K\"ahler metric on $B$. Given $\phi \in C^\infty(B)$, let $\omega_B = \omega_B' + i \delb\del \phi$. The aim is to find $\phi$ such that  $\omega_B$ solves equation (\ref{un-traced equation}). The reduction to the complex Monge--Ampere equation is standard; brief details are given below, for more information see, for example, the account in \cite{joyce:cmwsh}.

Since $\rho(\omega'_B) - a \in 2\pi(c_1(B) - \alpha)$, there exists $f$ such that 
$$
\rho(\omega'_B) - a = \lambda \omega'_B + i\delb\del f
$$
where $\lambda = 0$ if $a-c_1(B) = 0$ and $\lambda = -2\pi$ if $a -c_1(B) > 0$. Define $M(\phi) = \omega_B^r/\omega'^r_B$ where $r = \dim B$. Then
$
\rho(\omega_B) = \rho(\omega'_B) + i\delb\del \log M(\phi)
$.
Hence $\omega_B$ solves (\ref{un-traced equation}) if and only if $\phi$ solves
$
\log M(\phi) = \lambda \phi - f
$.
By the famous theorems of Aubin \cite{aubin:edtmaslvkc} ($\lambda = -1$) and Yau \cite{yau:otrcoackmatcmae1} ($\lambda = 0$), this has a unique solution. Hence there is a unique solution $\omega_B$ to (\ref{un-traced equation}) in $[\omega'_B]$.

It remains to show that $\omega_B$ is the unique solution in $[\omega'_B]$ to the \emph{traced} equation (\ref{traced equation}) appearing in Theorem \ref{old theorem}. Suppose $\omega_B$ solves $\Scal - \tr a = \lambda$, with $\lambda \omega_B \in a - c_1(B)$. Now $\rho(\omega_B) - a$ is $\delb$-closed and it follows from the K\"ahler identity $[\tr, \del] = i\delb^*$ and the fact that $\Scal -\tr a$ is constant that it is also $\delb^*$-closed. Hence $\rho(\omega_B) -a$ and $\lambda \omega_B$ are harmonic representatives for the same class and so, by Hodge theory, are equal. That is, $\omega_B$ solves (\ref{un-traced equation}) also. 

Theorem \ref{main theorem} now follows from Theorem \ref{old theorem}.
\end{proof}

\section{Non-negativity of the CM-line bundle}\label{non-negativity of alpha}

This section gives the proof of Theorem \ref{alpha is nef}. The key step is provided by the following theorem of Cornalba--Harris \cite{cornalba.harris:dcatfosv}.

\begin{theorem}[Cornalba--Harris]\label{cornalba-harris}
Let $\pi \colon X\to B$ be a flat proper morphism of varieties where $B$ is one-dimensional and $X$ is $(n+1)$-dimensional. Let $E \to X$ be a relatively very ample line bundle whose restriction to at least one fibre gives a Hilbert semi-stable embedding. 

Then 
$$
h c_1(E)^{n+1}-(n+1) c_1(E)^n \pi^* c_1(\pi_* E)
\geq 
0
$$
where $h$ is the rank of $\pi_*E$.
\end{theorem}

Theorem \ref{alpha is nef} is proved by looking at the asymptotics of this result.

\begin{proof}[Proof of Theorem \ref{alpha is nef}]

Recall that $\pi \colon X \to B$ is a holomorphic submersion between compact K\"ahler manifolds, $L$ is a relatively ample line bundle whose restriction to any fibre is asymptotically stable. The aim is to show that $\alpha$ (defined in (\ref{alpha})) is nef.
It suffices to consider the case when $B$ is a curve.

For sufficiently large $m$, $L^m$ is relatively very ample and makes at least one of the fibres Hilbert semi-stable. By Theorem \ref{cornalba-harris}, with $E = L^m$,
\begin{equation}\label{CH inequality}
h(m)m^{n+1}c_1(L)^{n+1} - m^n(n+1) c_1(L)^n \pi^*c_1(\pi_* (L^m))
\geq
0,
\end{equation}
where $h(m) = \rank \pi_*(L^m)$ is the Hilbert polynomial of the restriction of $L$ to a fibre.

By Grothendieck--Riemann--Roch, 
$$
c_1(\pi_* L^m)
=
\pi_*\left(
\frac{c_1(L)^{n+1}}{(n+1)!}\, m^{n+1}
-
\frac{c_1(L)^n c_1(K_{X/B})}{2 n!}\, m^n 
+
\cdots
\right).
$$
On the other hand,
$$
h(m)
=
\pi_*(c_1(L)^n) \left( \frac{1}{n!}m^n + \frac{s}{2n!}m^{n-1} + \cdots \right). 
$$
Apply $\pi_*$ to (\ref{CH inequality}) and collect terms. The $m^{2n+1}$-coefficient is zero and so the leading term is $m^{2n}$. The coefficient of $m^{2n}$ is $\alpha/(2(n+1)!)$. Since $\pi_*$ of (\ref{CH inequality}) is positive for all large $m$, its leading coefficient must be positive.
\end{proof}

One situation in which the positivity of $c_1(L_{\CM})$ is completely understood is that of curves of genus at least two. In this case, all polarisations are multiples of the canonical polarisation and $c_1(\CM)$ is essentially $\pi_*(c_1(K_{X/B})^2)$, the first tautological class, whose positivity properties are completely known. It is not ample on the Hilbert scheme, even over the stable locus (\cite{harris.morrison:moc}, page 313). On the other hand, it \emph{is} ample when considered over the smooth locus (\cite{harris.morrison:moc}, page 312). 

This tallies with the higher dimensional result of Fujiki--Schumacher \cite{fujiki.schumacher:tmsoeckmagwpm}, which shows that, provided the fibres have no holomorphic vector fields, $c_1(L_{\CM})$ is ample over the part of the smooth locus which admits a cscK metric. With this in mind, it seems reasonable to guess that $c_1(L_{\CM})$ is ample on the smooth, asymptotically semi-stable locus of the Hilbert scheme. 

Another question is whether or not the stability hypothesis in Theorem \ref{alpha is nef} is necessary. Certainly \emph{some} stability hypothesis is required. In \cite{cornalba.harris:dcatfosv} Cornalba and Harris discuss a family of Hilbert unstable varieties for which their inequality is violated. Is it possible, however, to replace asymptotic Hilbert semi-stability with K-semi-stability?

\section{Examples}\label{examples} 

A large class of submersions satisfying the hypotheses of Conjecture \ref{anti-fano base} are provided by the iterated surface bundles of Morita. Their construction is given in detail in chapter 4 of \cite{morita:gocc}. 

In what follows, a \emph{$\Sigma_g$-bundle} is a holomorphic submersion between complex manifolds whose fibres are curves of genus $g$. For each positive integer $k$, define the set $\mathcal C_k$ of connected compact complex $k$-folds recursively as follows: the only element of $\mathcal C_0$ is a single point; in general, $\mathcal C_{k+1}$ is defined to be the set consisting of any finite covering of the total space of a $\Sigma_g$-bundle with $g\geq 2$ and whose base belongs to $C_k$. Let $\mathcal C = \bigcup \mathcal C_k$ and call its members \emph{iterated surface bundles}. 

Given a $\Sigma_g$-bundle $\pi \colon X \to B$ with $g\geq2$ whose total space is an iterated surface bundle, Morita produces a finite covering $B' \to X$, and a $\Sigma_{g'}$-bundle $X' \to B'$ where $g' = m^2g - \frac{1}{2}m(m+1)+1$ is again at least 2. Moreover, the moduli of the fibres of $X'$ are not constant. (In \cite{morita:gocc}, Morita shows that certain characteristic classes of surface bundles are non-zero on $X'$; these classes vanish on isotrivial families.) It remains to check that these fibrations satisfy the hypotheses of Corollary \ref{anti-fano base}.

If $\pi \colon X \to B$ is a non-constant family of smooth curves of genus at least 2, then the relative canonical bundle $K_{X/B}$ is ample (see \cite{harris.morrison:moc}, page 309). If the canonical bundle of $B$ is also ample then so is the canonical bundle of $X$. That is, if $c_1(B) <0$ then $c_1(X) <0$. 

The short exact sequence
$
0\to V \to TX \to \pi^* TB \to 0
$,
where $V$ denotes the veritcal tangent bundle, gives a long exact sequence in cohomology
$
0 \to H^0(X, V) \to H^0 (X, TX)\to H^0(X, \pi^*TB) \to \cdots
$.
Now $H^0(X, V)=0$ as the fibres admit no non-zero holomorphic vector fields, whilst $H^0(X, \pi^*TB)=H^0(B, \pi_*\pi^*TB) = H^0(B,TB)$, so $H^0(B,TB)=0$ implies that $H^0(X, TX) = 0$. Similarly, finite covers of $B$ have no holomorphic vector fields.

Therefore repeatedly applying Morita's construction for various choices of $m$ and $g$ to $\Sigma_g$-bundles over a point produces infinitely many K\"ahler manifolds $X$ of arbitrary dimension all of which satisfy the hypotheses of Corollary \ref{anti-fano base}.

It should be pointed out that there are other ways of producing infinitely many cscK manifolds in each dimension, the most obvious being products. K\"ahler--Einstein metrics give non-trivial examples. For examples which are cscK but not K\"ahler--Einstein, the K\"ahler--Einstein metrics can be deformed via a theorem of LeBrun and Simanca \cite{lebrun.simanca:ekmacdt} which says that, in the absence of holomorphic vector fields, if a K\"ahler class admits a cscK representative, then nearby K\"ahler classes also admit a cscK representative. Alternatively, a result of Arezzo and Pacard \cite{arezzo.pacard:buadkoocsc} states that, in the absence of holomorphic vector fields, a cscK metric on $X$ gives a cscK metric on the blow-up of $X$ at a point, providing yet more examples. All these approaches, however, ultimately rely on deforming K\"ahler--Einstein metrics. By contrast, the examples produced via Theorem \ref{main theorem} lie in classes arbitrarily far from the canonical class and are unrelated to any K\"ahler--Einstein metric on $X$ or $B$.

\bibliographystyle{plain}
\bibliography{high_dim_fibs_bib}

\end{document}